\documentclass[12pt]{amsart}
\usepackage{amssymb}
\usepackage{amsmath}
\usepackage{amsthm}

\DeclareMathOperator{\Ad}{Ad}
\DeclareMathOperator{\ad}{ad}
\DeclareMathOperator{\dv}{div}
\newcommand{\A}{{\mathcal A}}

\DeclareMathOperator{\End}{End}

\renewcommand{\l}{{\mathcal L}}
\renewcommand{\r}{{\mathcal R}}
\newcommand{\N}{{\mathcal N}}

\newcommand{\R}{\ensuremath{\mathbb{R}}}

\newcommand{\p}{\partial}
\newcommand{\s}{{\rm Symb}}
\renewcommand{\S}{\Sigma}
\newcommand{\T}{{T^\ast M}}
\newcommand{\D}{{\mathcal D}}

\newtheorem{lemma}{Lemma}
\newtheorem{proposition}{Proposition}
\newtheorem{theorem}{Theorem}
\newtheorem{corollary}{Corollary}
\begin{document}
\title[Inverse mapping of formal symplectic groupoid]{On the inverse mapping of the formal symplectic groupoid of a deformation quantization}
\author[A.V. Karabegov]{Alexander V. Karabegov}\footnote{Research was partially supported by an ACU Math/Science grant.}
\address[Alexander V. Karabegov]{Department of Mathematics and Computer Science, Abilene Christian University, ACU Box 28012, 252 Foster Science Building, Abilene, TX 79699-8012}
\email{alexander.karabegov@math.acu.edu}

\begin{abstract}
To each natural star product on a Poisson manifold $M$ we associate an antisymplectic involutive automorphism of the formal neighborhood of the zero section of the cotangent bundle of $M$. If $M$ is symplectic, this mapping is shown to be the inverse mapping of the formal symplectic groupoid of the star product. The construction of the inverse mapping involves modular automorphisms of the star product.
\end{abstract}

\date{February 27, 2004}
\maketitle

\section{Introduction}

Symplectic groupoids are geometric objects of semiclassical origin. They were introduced independently by Karas\"ev \cite{Ka}, Weinstein \cite{W} and Zakrzewski \cite{Z}. Recently symplectic groupoids appeared in the context of deformation quantization in \cite{CF2},\cite{CF3},\cite{CDF}. In \cite{CDF} A. S. Cattaneo, B. Dherin, and G. Felder considered the formal integration problem for Poisson manifolds whose solution is given by a formal symplectic groupoid. They start with the zero Poisson structure on $\R^n$, the corresponding trivial symplectic groupoid $T^* \R^n$, and a generating function of the Lagrangian product space of this groupoid. A formal symplectic groupoid is then defined in terms of a formal deformation of that trivial generating function.
One of the main results of  \cite{CDF} is an explicit formula for the generating function that defines the formal symplectic groupoid related to Kontsevich star product.

One can also take an equivalent approach to the definition of a formal symplectic groupoid. If $\Sigma$ is a symplectic groupoid with the unit space $\Gamma \subset \Sigma$, one can consider the formal neighborhood $(\Sigma,\Gamma)$ of $\Gamma$ in $\Sigma$ as a formal symplectic groupoid (the definition of a formal neighborhood is given in Section \ref{S:natural}, see also the Appendix to \cite{Deq}). The groupoid operations and axioms have then their natural formal counterparts that can be used to state a definition of a formal symplectic groupoid.

The simlest example of a symplectic groupoid is the pair symplectic groupoid $M \times \bar M$ of a symplectic manifold $M$. Here $\bar M$ is a copy of $M$ with the opposite symplectic structure. In \cite{Deq} we related to any natural deformation quantization on a symplectic manifold $M$ a canonical formal symplectic groupoid isomorphic to the formal pair symplectic groupoid $(M \times \bar M, M_{\rm diag})$, where $M_{\rm diag}$ is the diagonal of $M \times \bar M$. The isomorphism was given via the source and target mappings of the formal symplectic groupoid of the deformation quantization. The construction of these mappings remains valid in the Poisson case.
The results of \cite{CDF} strongly suggest that one can relate a canonical formal symplectic groupoid to any natural star product on a Poisson manifold.

In this Letter we give a construction of the inverse mapping of the formal symplectic groupoid associated to a natural deformation quantization on a symplectic manifold that still makes sense in the Poisson case. The construction uses modular automorphisms of the deformation quantization in question. Modular automorphisms in deformation quantization were introduced in \cite{Wald} in the framework of the Tomita-Takesaki theory where  the construction of the automorphisms involves complex conjugation.
Our construction is rather in the spirit of  \cite{Wmod}. In particular, the coefficient at the first power of the formal parameter  of a modular automorphism is a modular vector field. We also show that the modular automorphisms of a given star product determine a single outer automorphism of the star product which is an obstruction to the existence of a nonvanishing trace density of that star product.

\section{Transpositions of formal differential operators}\label{S:trans}

Denote by $\D$ the algebra of differential operators on a real manifold $M$ and by $\D[[\nu]]$ the algebra of formal differential operators on $M$. Let $\rho$ be a nonvanishing smooth density on $M$. 
It determines the bilinear pairing 
\begin{equation}\label{E:pair}
   (u, v)_\rho = \int_M u \cdot v\, \rho
\end{equation}
of $C^\infty(M)$ with the space of compactly supported smooth functions $C_c^\infty(M)$. For a differential operator $A$ on $M$ one can define its formal transpose $A^t_\rho \in \D$ with respect to that pairing so that $(Au,v)_\rho = (u, A^t_\rho v)_\rho$ for all $u \in C^\infty(M)$ and $v \in C_c^\infty(M)$. We will call a formal density $\rho = \rho_0 + \nu \rho_1 + \ldots$ on $M$ nonvanishing if $\rho_0$ is nonvanishing. For a nonvanishing formal density $\rho$ on $M$ one can define a $\nu$-linear pairing $(\cdot,\cdot)_\rho$ of $C^\infty(M)[[\nu]]$ and $C^\infty_c(M)[[\nu]]$ by the same formula (\ref{E:pair}) and the formal transpose $A^t_\rho$ of a formal differential operator $A\in \D[[\nu]]$. If $\tilde \rho = \phi \rho$ is another nonvanishing formal density on $M$ (so that $\phi$ is an invertible formal function on $M$) then for $u \in C^\infty(M)[[\nu]], \ v \in C^\infty_c(M)[[\nu]]$, and $A \in \D[[\nu]]$ we have $\int Au \cdot v \, \tilde\rho = \int Au \cdot v \phi\, \rho = \int u \cdot A^t_\rho (\phi v)\, \rho$. On the other hand, $\int Au \cdot v \, \tilde\rho =  \int u \cdot A^t_{\tilde\rho}v \, \tilde\rho = \int  u \cdot \phi A^t_{\tilde\rho}v \, \rho$. Therefore
\begin{equation}\label{E:tilderho}
         A^t_{\tilde\rho} = \phi^{-1} A^t_\rho \phi.
\end{equation}

Deformation quantization on a Poisson manifold $(M, \{\cdot,\cdot\})$, as introduced in \cite{BFFLS},  is an associative algebra structure on the space of formal functions $C^\infty(M)[[\nu]]$. Its product (named a star product) is given on $f,g\in
C^\infty(M)$ by the formula
\begin{equation}\label{E:star}
f \ast g = \sum_{r = 0}^\infty \nu^r C_r(f,g), 
\end{equation} 
where $C_r,\, r\geq 0,$ are bidifferential operators on $M$, $C_0(f,g) =
fg$ and $C_1(f,g) - C_1(g,f) = \{f,g\}$.  We adopt the convention that
the unit of a star-product is the unit constant. Two differentiable
star-products $\ast,\ast'$ on a Poisson manifold $(M,\{\cdot,\cdot\})$ are
called equivalent if there exists an isomorphism of algebras $B:
(C^\infty(M)[[\nu]],\ast) \to (C^\infty(M)[[\nu]],\ast')$ of the form $B=
1 + \nu B_1 + \nu^2 B_2 + \dots,$ where $B_r, r\geq 1,$ are differential
operators on $M$. The existence and classification problem for deformation quantization was first solved in the non-degenerate (symplectic) case (see \cite{DWL}, \cite{OMY}, \cite{F1} for existence proofs and \cite{F2}, \cite{NT}, \cite{D}, \cite{BCG}, \cite{X} for classification) and then Kontsevich \cite{K} showed that every Poisson manifold admits a deformation quantization and that the equivalence classes of deformation quantizations can be parameterized by the formal deformations of the Poisson structure. 

Let $\A = (C^\infty(M)[[\nu]],*)$ be a deformation quantization on a Poisson manifold $(M,\{\cdot,\cdot\})$ and  $\l ,\r \subset \D[[\nu]]$ be the algebras of operators of $*$-multiplication from the left and from the right, respectively. For $f,g \in C^\infty(M)[[\nu]]$ denote by $L_f \in \l$ and $R_g\in \r$ the operators of $*$-multiplication from the left by $f$ and from the right by $g$, respectively, so that $L_f g = f \ast g = R_gf$. The associativity of $\ast$ is equivalent to the fact that $[L_f, R_g]=0$ for any formal functions $f,g$ on $M$. It is well known that $\l$ and $\r$ are commutants of each other in $\D[[\nu]]$. 
Introduce an operator $J_\rho \in \D[[\nu]]$ such that $J_\rho f = (R_f)^t_\rho 1$. If $f\in C^\infty(M)$ then $R_f = f \pmod{\nu}$, whence $J_\rho = 1 \pmod{\nu}$. Therefore $J_\rho$ is invertible and $J_\rho^{-1} = 1 \pmod{\nu}$.   For $u \in C^\infty(M)$ and $v \in C_c^\infty(M)$ we get
\begin{equation}\label{E:trans}
    \int_M u * v \, \rho = \int_M R_v u\, \rho = \int_M u \cdot (R_v)^t_\rho 1\, \rho = \int_M u \cdot J_\rho v \, \rho.
\end{equation}
Define an operator $Q_\rho \in \D[[\nu]]$ by the formula $Q_\rho = (J_\rho^{-1})^t_\rho J_\rho$. Using formula (\ref{E:trans}) we obtain that \begin{equation}\label{E:qtrace}
      \int_M Q_\rho u * v \, \rho = \int_M  (J_\rho^{-1})^t_\rho J_\rho u \cdot J_\rho v \, \rho = \int_M v \cdot J_\rho u \, \rho = \int_M v * u \, \rho.
\end{equation}
It follows from (\ref{E:trans}) and (\ref{E:qtrace}) that $\rho$ is a trace density of the star product $*$ if and only if $J_\rho$ is formally symmetric, i.e., $(J_\rho)^t_\rho = J_\rho$, or, equivalently, $Q_\rho = (J_\rho^{-1})^t_\rho J_\rho = 1$. Assume that $\tilde\rho$ is another formal smooth trace density  of $*$ that is not necessarily nonvanishing. Then $\tilde\rho = \phi \rho$ for some formal function $\phi \in C^\infty(M)[[\nu]]$. It turns out that all such trace densities are naturally parameterized by the elements of the center of the algebra $\A$ (the Casimir elements). For any $\tilde\rho = \phi \rho, \ u \in C^\infty(M)$, and $v \in C_c^\infty(M)$ we have from (\ref{E:trans}) that  $\int u * v\ \tilde\rho = \int u * v \cdot \phi\, \rho = \int u * v * J_\rho^{-1}\phi \ \rho = \int u \cdot J_\rho (v * J_\rho^{-1}\phi) \,\rho$. On the other hand, using the fact that $\rho$ is a trace density, we get $\int v * u\ \tilde\rho = \int v * u \cdot \phi\, \rho = \int v * u * J_\rho^{-1}\phi \ \rho = \int  u * J_\rho^{-1}\phi * v \ \rho = \int u \cdot J_\rho (J_\rho^{-1}\phi * v) \,\rho$. Thus $\tilde\rho$ is a trace density of the star product $*$ if and only if $v * J_\rho^{-1}\phi = J_\rho^{-1}\phi * v$ for any $v \in C_c^\infty(M)$, i.e., iff $\psi = J_\rho^{-1}\phi$ is a Casimir function of $*$. Therefore the trace densities of the star product $*$ are of the form $(J_\rho \psi) \, \rho$, where $\psi$ is a Casimir function of $*$. 

Introduce a mapping $K_\rho \in \End \D[[\nu]]$ by the formula
$K_\rho [A] = J_\rho^{-1} A^t_\rho J_\rho$, where $A \in \D[[\nu]]$.
Then for any $u \in C^\infty(M),\ v \in C^\infty_c(M)$ and $A \in \D[[\nu]]$ we have, using (\ref{E:trans}), that $\int Au * v \, \rho = \int Au \cdot J_\rho v \, \rho = \int u \cdot A^t_\rho J_\rho v \, \rho = \int u *J^{-1}_\rho A^t_\rho J_\rho v \, \rho = \int u * K_\rho[A]v\, \rho$. Thus
\begin{equation}\label{E:ktransp}
            \int Au * v \, \rho = \int u * K_\rho[A]v\, \rho,
\end{equation}
i.e., $K_\rho[A]$ is the transposition of $A$ with respect to the pairing $(u,v) \mapsto \int u * v \, \rho$. 
The square of the mapping $K_\rho$ is the conjugation by the operator $Q_\rho$,
\begin{equation}\label{E:square}
  \ K_\rho\left[K_\rho[A]\right] = J_\rho^{-1}\left(J_\rho^{-1} A^t_\rho J_\rho\right)^t_\rho J_\rho = J_\rho^{-1}(J_\rho)^t_\rho A (J^{-1}_\rho)^t_\rho J_\rho = Q_\rho^{-1} A Q_\rho.
\end{equation}

We need to show that $Q_\rho 1 = 1$. It is easy to check that $J_\rho1 = 1$. Now it is sufficient to show that $(J_\rho)^t_\rho 1 = 1$ as well. Using formula (\ref{E:trans}) we obtain for $v \in C^\infty_c(M)$ that $\int v \, \rho = \int 1 * v \, \rho = \int  1 \cdot J_\rho v \, \rho = \int  (J_\rho)^t_\rho 1 \cdot v \, \rho$, whence the statement follows.

Since for $f,u \in C^\infty(M)[[\nu]]$ and $v \in C^\infty_c(M)[[\nu]]$ we have $\int R_f u * v \, \rho = \int u * f * v \, \rho = \int u * L_f v \, \rho$, it follows that 
\begin{equation}\label{E:krl}
K_\rho[R_f] = L_f. 
\end{equation}
Thus $J_\rho^{-1}(R_f)^t_\rho J_\rho = L_f$, or, equivalently, 
\begin{equation}\label{E:jlr}
                            J_\rho L_f J_\rho^{-1} = (R_f)^t_\rho.
\end{equation}
For $g\in C^\infty(M)[[\nu]]$ we obtain from (\ref{E:jlr}) that $(L_g)^t_\rho =(J_\rho)^t_\rho R_g (J_\rho^{-1})^t_\rho$. Since $R_f$ commutes with $L_g$ we get that $J_\rho L_f J_\rho^{-1}$ commutes with $(J_\rho)^t_\rho R_g (J_\rho^{-1})^t_\rho$, which implies that 
$A = (J_\rho^{-1})^t_\rho J_\rho L_f J_\rho^{-1} (J_\rho)^t_\rho = Q_\rho L_f Q_\rho^{-1}$ commutes with $R_g$ for arbitrary $g$. This means that $A \in \l$, i.e.,  that $\l$ is invariant with respect to the conjugation by the operator $Q_\rho$. Thus the conjugation by the operator $Q_\rho$ establishes an automorphism of the algebra $\l$. Since the mapping $f \mapsto L_f$ is an isomorphism of  the algebra $\A$ onto $\l$, we obtain the following automorphism of the star product $*$, $f \mapsto Q_\rho L_f Q_\rho^{-1}1 = Q_\rho f$. Here we have used that $Q_\rho 1 = 1$. Thus $Q_\rho$ is an automorphism of the star product $*$. This is  the modular automorphism of the deformation quantization $\A$ corresponding to the formal density $\rho = \rho_0 + \nu \rho_1 + \ldots$. The following calculation shows that the derivation $\log Q_\rho$ of the algebra $\A$ is a deformation of the modular vector field $C^\infty(M) \ni  f \mapsto \dv_{\rho_0} H_f$ on the Poisson manifold $(M,\{\cdot,\cdot\})$ corresponding to the density $\rho_0$, introduced in \cite{Wmod}. Here $H_f$ denotes the Hamiltonian vector field corresponding to $f$ with respect to the Poisson structure on $M$, so that $H_f g = \{f,g\}$.
For $f \in C^\infty(M)$ set $g = (J_\rho^{-1})^t_\rho J_\rho f = Q_\rho f$. Thus $J_\rho f = (J_\rho)^t_\rho g$ and the function $g$ can be represented as $g = f + \nu h \pmod{\nu^2}$.  For $u \in C_c^\infty(M)$ we have $\int u * f \, \rho = \int u \cdot J_\rho f \, \rho = \int u \cdot (J_\rho)^t_\rho g \, \rho = \int g \cdot J_\rho u \, \rho =
\int g * u \, \rho$. Now
\begin{equation}\label{E:lft}
\int u * f \, \rho = \int uf \, \rho_0 + \nu\left(\int C_1(u,f)\, \rho_0 + \int uf \, \rho_1\right)  \pmod{\nu^2}. 
\end{equation}
Similarly,
\begin{equation}\label{E:rght}
\int g * u \, \rho = \int fu \, \rho_0 + \nu\left(\int C_1(f,u) + hu\, \rho_0 + \int fu \, \rho_1 \right)  \pmod{\nu^2}. 
\end{equation}
Subtracting (\ref{E:lft}) from (\ref{E:rght}) and using that  $C_1(f,u) - C_1(u,f) = \{f ,u\} = H_f u$  we get that $\int H_f u + hu \, \rho_0 = 0$. Since $\int H_f u \, \rho = - \int u \dv_{\rho_0} H_f \, \rho_0$, we finally obtain that $h = \dv_{\rho_0} H_f$ which implies that $\log Q_\rho f = \nu \dv_{\rho_0} H_f \pmod{\nu^2}$.

Let $\rho$ and $\tilde\rho = \phi\, \rho$ be nonvanishing formal densities on $M$. Denote $\psi = J_\rho^{-1}\phi$. For $u \in C^\infty(M)$ and $v \in C^\infty_c(M)$ we have, using formula (\ref{E:trans}), that  $\int u * v \, \tilde\rho = \int u * v \cdot \phi \, \rho = \int u * v * \psi \, \rho = 
\int u \cdot J_\rho(v * \psi) \, \rho = \int u \cdot J_\rho R_\psi \, \rho$. On the other hand, $\int u * v \, \tilde\rho = \int u \cdot J_{\tilde\rho} v \, \tilde\rho = \int u \cdot \phi J_{\tilde\rho} v \, \rho$, whence $\phi J_{\tilde\rho} = J_\rho R_\psi$ or, equivalently,
\begin{equation}\label{E:jtilderho}
     J_{\tilde\rho} = \phi^{-1} J_\rho R_\psi.
\end{equation}
Using formulas (\ref{E:tilderho}) and (\ref{E:jtilderho}), we can express $K_{\tilde\rho}$ via $K_\rho$. For $A \in \D[[\nu]]$ we have that
$ K_{\tilde\rho}[A] = J^{-1}_{\tilde\rho} A^t_{\tilde\rho} J_{\tilde\rho} = (\phi^{-1} J_\rho R_\psi)^{-1} \phi^{-1} A^t_\rho \phi (\phi^{-1} J_\rho R_\psi) =(R_\psi)^{-1} J_\rho^{-1} A^t_\rho J_\rho R_\psi =   (R_\psi)^{-1}K_\rho[A] R_\psi.$ Thus
\begin{equation}\label{E:ktilderho}
    K_{\tilde\rho}[A] = (R_\psi)^{-1}K_\rho[A] R_\psi.
\end{equation}
Finally, we express $Q_{\tilde\rho}$ via $Q_\rho$. From (\ref{E:jtilderho}) we have that $(J_{\tilde\rho}^{-1})^t_{\tilde\rho} = \phi \circ (J_\rho^{-1})^t_{\tilde\rho} (R_\psi^{-1})^t_{\tilde\rho}$. Now, using formulas (\ref{E:tilderho}) and (\ref{E:jlr}), we get 
\begin{eqnarray*}
\lefteqn{Q_{\tilde\rho} =  (J_{\tilde\rho}^{-1})^t_{\tilde\rho}J_{\tilde\rho} =
\phi \circ (J_\rho^{-1})^t_{\tilde\rho} (R_\psi^{-1})^t_{\tilde\rho} \circ 
\phi^{-1}\circ J_\rho R_\psi = }\\
&& (J_\rho^{-1})^t_\rho (R_\psi^{-1})^t_\rho J_\rho R_\psi =
(J_\rho^{-1})^t_\rho J_\rho (L_\psi)^{-1} J_\rho^{-1} J_\rho R_\psi =
Q_\rho (L_\psi)^{-1} R_\psi.
\end{eqnarray*} Therefore,
\begin{equation}\label{E:qtilderho}
     Q_\rho = Q_{\tilde\rho}(R_\psi)^{-1}L_\psi  = Q_{\tilde\rho} \circ \Ad(\psi).
\end{equation}
Formula (\ref{E:qtilderho}) shows that the modular automorphisms $Q_\rho$ of the star algebra $\A$ corresponding to different nonvanishing densities $\rho$ differ by an inner automorphism of $\A$ and thus define a canonical outer automorphism of $\A$. This outer automorphism is an obstruction to the existence of a nonvanishing trace density. It was observed above that $\rho$ is a trace density iff $Q_\rho = 1$. Then, according to (\ref{E:qtilderho}), all other modular automorphisms of $\A$ are inner. Now assume that the modular automorphism  $Q_\rho$ is inner, say, $Q_\rho = \Ad (\chi)$. Then formula (\ref{E:qtilderho}) implies that  if $\psi = \chi$ then $Q_{\tilde\rho} =1$ and therefore $\tilde\rho = (J_\rho\chi) \rho$ is a trace density.

For a given star product $*$ and a nonvanishing formal density $\rho$ on the Poisson manifold $(M,\{\cdot,\cdot\})$ the operator $J_\rho$ can be used as an equivalence operator with a  star product $\tilde *$ on $(M,\{\cdot,\cdot\})$ given by the formula $f \tilde * g = J_\rho (J_\rho^{-1}f \ast J_\rho^{-1}g)$. Denote by $\tilde L_f$ and $\tilde R_f$ the operators of multiplication from the left and from the right by $f \in C^\infty(M)[[\nu]]$ with respect to the star product $\tilde *$, respectively. Using formula (\ref{E:trans}) we obtain for any $f,g \in C^\infty(M)[[\nu]]$ and $u \in C_c^\infty(M)$ the following chain of equalities, $\int u \cdot \tilde L_f g \, \rho = \int u \cdot (f \tilde * g) \, \rho = \int u \cdot J_\rho(J_\rho^{-1}f * J_\rho^{-1}g) \, \rho = \int u * (J_\rho^{-1}f * J_\rho^{-1}g) \, \rho = \int (u * J_\rho^{-1}f) * J_\rho^{-1}g \, \rho = \int (u * J_\rho^{-1}f) \cdot g \, \rho = \int R_{J_\rho^{-1}f} u \cdot g \, \rho = \int  u \cdot (R_{J_\rho^{-1}f})^t_\rho g \, \rho$, whence 
\begin{equation}\label{E:ltransr}
\tilde L_f = (R_{J_\rho^{-1}f})^t_\rho. 
\end{equation}

\section{Natural star products}\label{S:natural}

Let $M$ be a real $m$-dimensional manifold. Its cotangent manifold $\T$ is endowed with the standard Poisson bracket $\{\cdot,\cdot\}_\T$ 
given in local coordinates  by the formula
\[
\{f,g\}_\T = \frac{\p f}{\p \xi_k}\frac{\p g}{\p x^k} -  \frac{\p g}{\p \xi_k}\frac{\p f}{\p x^k}.
\]
Here as usual $\{\xi_l\}$ are the fibre coordinates dual to the base coordinates $\{x^k\}$. The algebra  $\D$ of differential operators on $M$
has a natural filtration with respect to the order of differential operators, $0 \subset \D_1 \subset \D_2 \subset \ldots$.
The symbol mapping $\s_r: \D_r \to C^\infty(\T)$ maps the operators of order not greater than $r$ to the functions on $\T$ that are homogeneous polynomials of degree $r$ on the fibres. The kernel of $\s_r$ is $\D_{r-1}$. In local  coordinates an operator $X\in\D_r$ with the leading term $X^J \p_J$ where $J$ runs over the multi-indices of order $r$ and $\p_k = \p/ \p x^k$, has the symbol $\s_r (X) = X^J \xi_J$. 

 For two arbitrary differential operators $X\in\D_k$ and $Y\in \D_l$ it is well known that $XY \in \D_{k+l},\ [X,Y] \in \D_{k+l-1}$ and
\begin{equation}\label{E:product}
\s_{k+l} XY = \s_k X\cdot \s_l Y,
\end{equation}
\begin{equation}\label{E:comm} 
\s_{k+l-1}([X,Y]) = \{\s_k(X),\s_l(Y)\}_\T.    
\end{equation}

A formal differential operator $A = A_0 + \nu A_1 + \nu^2 A_2+ \ldots \in \D[[\nu]]$ such that the order of $A_r$ is not greater than $r$ will be called {\it natural}. Denote by $\N$ the algebra of natural formal differential operators on $M$. Formula (\ref{E:comm}) implies that for any $A, B \in \N$ the operator $\frac{1}{\nu}[A,B]$ is natural. Denote by ${\mathfrak g}$ the pronilpotent Lie algebra of  the formal differential operators on $M$ that can be expressed as $\frac{1}{\nu}A$ where $A = \nu^2 A_2 + \nu^3 A_3 +\dots$, i.e., $A = 0 \pmod{\nu^2}$.

Let $Y$ be a closed $k$-dimensional submanifold of a real $n$-dimensional manifold $X$ and $I_Y \subset C^\infty(M)$ be the ideal of smooth functions on $X$ vanishing on $Y$. Then the quotient algebra $C^\infty(X,Y) := C^\infty(X)/I_Y^\infty$ can be thought of as the algebra of smooth functions on the formal neighborhood $(X,Y)$ of the submanifold $Y$ in $X$. If $U \subset X$ is a local coordinate chart with coordinates $\{x^i\}$ such that $U\cap Y$ is given by the equations $x^{k+1} = 0,\ldots,x^n =  0$, then $C^\infty(U,U \cap Y)$ is isomorphic to $C^\infty(U\cap Y)[[x^{k+1},\ldots, x^n]]$, where the isomorphism is established via the formal Taylor expansion of the functions on $U$ in the variables $x^{k+1},\ldots, x^n$.

Define the $\sigma$-symbol of an element $A = A_0 + \nu A_1 + \nu^2 A_2+ \ldots\in \N$ by the formula $\sigma(A) = \s_0 (A_0) + \s_1(A_1) + \ldots$. Such a formal series of homogeneous functions on $\T$ can be treated as an element of the algebra $C^\infty(\T,Z)$, where $Z$ is the zero section of $\T$.  The Poisson bracket $\{\cdot,\cdot\}_\T$ can be transferred to $C^\infty(\T,Z)$. We will use the same notation $\{\cdot,\cdot\}_\T$ for the bracket on $C^\infty(\T,Z)$.
Since $Z$ is a Lagrangian manifold in $\T$, the ideal $I_Z \subset C^\infty(\T)$ is closed with respect to the Poisson bracket $\{\cdot,\cdot\}_\T$. Set $I_Z^0 = C^\infty(\T)$. It is easy to check that for any $n \geq 0$ we have the inclusions $\{I_Z^2, I_Z^n\}_\T \subset I_Z^{n+1}$. Therefore the formal functions from $C^\infty(\T,Z)$ that vanish to second order on the zero section $Z$ form a pronilpotent Lie algebra ${\mathfrak h} := I_Z^2/I_Z^\infty$ with respect to the Poisson bracket $\{\cdot,\cdot\}_\T$. The formal Hamiltonian vector fields of the elements of  ${\mathfrak h}$ can be exponentiated to formal symplectomorphisms of $C^\infty(\T,Z)$.

Introduce the $\sigma$-symbol mapping $\sigma: \N \to C^\infty(\T,Z)$ that maps natural operators to their $\sigma$-symbols. Notice that the kernel of the mapping $\sigma$ is $\nu\N$. According to formulas (\ref{E:product}) and (\ref{E:comm}), the mapping $\sigma: \N \to C^\infty(\T,Z)$ is an algebra  homomorphism such that for any $A, B \in \N$
\begin{equation}\label{E:commab}
     \sigma\left(\frac{1}{\nu}[A,B]\right) = \{\sigma(A), \sigma(B)\}_\T = H_{\sigma(A)}\sigma(B),
\end{equation}
where $H_{\sigma(A)}$ is the (formal) Hamiltonian vector field corresponding to the $\sigma$-symbol of the operator $A$.

For a natural operator $X$ such that $X = 0 \pmod{\nu^2}$ its $\sigma$-symbol $\sigma(X)$ vanishes to second order on the zero section $Z$ of $\T$. According to formula (\ref{E:commab}), the mapping ${\mathfrak g} \ni \frac{1}{\nu} X \mapsto \sigma(X)$ defines a Lie algebra homomorphism of ${\mathfrak g}$ to ${\mathfrak h}$.

Consider the action of an invertible formal differential operator $B \in \D[[\nu]]$ on $\D[[\nu]]$ by conjugation, $\D[[\nu]] \ni A \mapsto BAB^{-1}$.
\begin{lemma}\label{L:conj} 
\noindent (i) If $B$ is a natural formal differential operator, its action by conjugation leaves the algebra $\N$ of natural operators invariant and it induces the trivial action on the $\sigma$-symbols.

\noindent (ii) If $B$ can be represented in the form $B = \exp \frac{1}{\nu} X$, where $X$ is a natural differential operator such that $X = 0 \pmod{\nu^2}$, then the action of the operator $B$ by conjugation leaves the algebra $\N$ of natural operators invariant and it induces the action on the $\sigma$-symbols by the operator $\exp H_{\sigma (X)}$, where $H_{\sigma (X)}$  is the (formal) Hamiltonian vector field corresponding to the $\sigma$-symbol of the operator $X$.
\end{lemma}
\begin{proof} (i) Since $\sigma(1) = 1$, we see that $\sigma(B^{-1}) = \left(\sigma (B)\right)^{-1}$. For any natural operator $A \in \N$ the operator  $BAB^{-1}$ is also natural and
$\sigma(BAB^{-1}) = \sigma(B) \sigma (A) 
\left(\sigma (B)\right)^{-1} = \sigma(A)$.

\noindent (ii) Since $X = 0 \pmod{\nu^2}$, for any formal differential operator $A \in \D[[\nu]]$ the series 
\[
BAB^{-1} = \sum_{n = 0}^\infty \frac{1}{n!}\left( \frac{1}{\nu} 
\ad(X) \right)^n A
\]
is $\nu$-adically convergent. It follows from formula (\ref{E:comm}) that
if the operator $A$ is natural, then the operator $\left( \frac{1}{\nu} \ad(X) \right)^n A$ is natural as well. Formula
(\ref{E:commab}) implies that
\[ 
\sigma\left(\left( \frac{1}{\nu} \ad(X) \right)^n A\right) = \left(H_{\sigma(X)}\right)^n \sigma (A), 
\]
whence
\begin{equation}\label{E:sigmaconj}
   \sigma(BAB^{-1}) = \left(\exp H_{\sigma(X)} \right)\sigma (A).
\end{equation}
\end{proof}

It is well known that all the explicit constructions of star-products enjoy the following property: for all $r \geq 0$ the bidifferential operator $C_r$ in (\ref{E:star}) is of order not greater than $r$ in both arguments (most important examples are Fedosov's star-products on symplectic manifolds and  Kontsevich's star-product on $\R^n$ endowed with an arbitrary Poisson bracket). The star-products with this property were called natural by Gutt and Rawnsley in \cite{GR}, where general properties of such star-products were studied (before these star-products were said to be of Vey type). 

A star-product $\ast$ on $M$ is natural if and only if for any  $f,g \in C^\infty(M)[[\nu]]$ the operators $L_f, R_g$ are natural. It turns out that for a star-product $\ast$ to be natural it is sufficient that $L_f \in \N$  for all $f \in C^\infty(M)[[\nu]]$ or, similarly, that $R_f \in \N$  for all $f \in C^\infty(M)[[\nu]]$. 
\begin{proposition}\label{P:nat}
If for any $f \in C^\infty(M)[[\nu]]$ the operator $L_f$ of a star-product $\ast$ is natural, then the star-product $\ast$ is also natural. Similarly, if for any $f \in C^\infty(M)[[\nu]]$ the operator $R_f$ is natural, then the star-product $\ast$ is also natural.
\end{proposition}
\begin{proof} We will prove the first statement of the proposition only. Fix an arbitrary formal function $g \in C^\infty[[\nu]]$. For any $f\in C^\infty(M)$ the operator $L_f  = f + \nu A_1 + \nu^2 A_2 +\ldots$ is natural  and commutes with $R_g = B_0 + \nu B_1 + \nu^2 B_2 + \ldots$. We claim that the operator $R_g$ is also natural. Since $[L_f,R_g] = 0$, we have for each $n \geq 1$ that
\begin{equation}\label{E:natq}
      [B_n , f] = \sum_{k =1} ^ n [A_k, B_{n-k}]
\end{equation}
and that $B_0$ commutes with the point-wise multiplication operator by an arbitrary function, whence it is a multiplication operator itself. Now assuming that the order of $B_r$ is not greater than $r$ for all $r < n$, we get from (\ref{E:natq}) that this is true for $r=n$ as well. The statement of the Proposition follows by induction.
\end{proof}

Let $*$ be a natural star product and  $\rho$ a nonvanishing formal density on the Poisson manifold $(M,\{\cdot,\cdot\})$. Define a  star product $\tilde *$ on $(M,\{\cdot,\cdot\})$ by the formula $f \tilde * g = J_\rho (J_\rho^{-1}f \ast J_\rho^{-1}g)$. We see from formula (\ref{E:ltransr}) that the operator $\tilde L_f$ is natural for any $f \in C^\infty(M)[[\nu]]$. Thus we obtain from Proposition \ref{P:nat} that the star product $\tilde *$ is natural as well.

In \cite{GR} the following important theorem was proved.
\begin{theorem}\label{T:GR} {\rm (S. Gutt, J. Rawnsley, \cite{GR})}
If an equivalence operator $B$ of two natural star products on a Poisson manifold $M$ is expressed in the form $B = \exp\frac{1}{\nu} X$, where $X$ is a formal differential operator on $M$ such that  $X = 0 \pmod{\nu^2}$, then $X$ is natural.
\end{theorem}
Theorem \ref{T:GR} implies that the equivalence operator $J_\rho$ of the natural star products $*$ and $\tilde *$ can be expressed in the form $J_\rho = \exp\frac{1}{\nu} X$, where the operator $X$ is natural and $X = 0 \pmod{\nu^2}$. Therefore the mapping $\N \ni A \mapsto J_\rho A J_\rho^{-1}$ is an automorphism of the algebra $\N$ of natural operators.
It induces  via the $\sigma$-mapping the automorphism of $C^\infty(\S,Z)$ by the operator $\exp H_{\sigma(X)}$, where $H_{\sigma(X)}$ is the (formal) Hamiltonian vector field corresponding to the $\sigma$-symbol of the operator $X$.  Now let $\tilde \rho$ be another nonvanishing formal density on $M$. Then $\tilde\rho = 
\phi \, \rho$ for some invertible formal function $\phi \in C^\infty(M)[[\nu]]$.
Formula \ref{E:jtilderho} and Lemma \ref{L:conj} imply that the conjugation via the operators $J_\rho$ and $J_{\tilde\rho}$ induces the same mapping on the $\sigma$-symbols, which means that the
automorphism $\exp H_{\sigma(X)}$ is independent of the choice of the density $\rho$.

In \cite{Deq} we associated to each natural star product $*$ on a Poisson manifold $(M,\{\cdot,\cdot\})$ the mappings $S,T:C^\infty(M) \to C^\infty(\T,Z)$ such that for $f\in C^\infty(M)\ Sf = \sigma(L_f)$ and $Tf = \sigma(R_f)$. It was shown that $S$ is a Poisson morphism, $T$ is an anti-Poisson morphism, and that for any $f,g\in C^\infty(M) \ Sf = f \pmod{\xi}, Tg = g \pmod{\xi}$, and $\{Sf,Tg\}_\T = 0$. 
The mappings $S$ and $T$ are induced by the formal Poisson morphism $s: (\T,Z) \to M$ and anti-Poisson morphism $t: (\T,Z) \to M$, respectively. Denote by $\bar M$ a copy of the manifold $M$ endowed with the opposite Poisson structure. Let  the manifold $M \times \bar M$ be endowed with the product Poisson structure and $M_{\rm diag}$ denote its diagonal. The mapping $S \otimes T:C^\infty(M \times \bar M, M_{\rm diag}) \to C^\infty(\T,Z)$ is therefore a Poisson morphism. If $M$ is symplectic, $S \otimes T$ is an isomorphism and thus the dual
morphism $s \times t : (\T,Z) \to (M \times \bar M, M_{\rm diag})$ is a formal symplectic isomorphism.
The formal manifold $(M \times \bar M, M_{\rm diag})$ is a formal (pair) symplectic groupoid. Hence $(\T,Z)$ can be thought of as an isomorphic formal symplectic groupoid for which $s$ and $t$ are the source and target morphisms, respectively. The inverse mapping of the pair groupoid $M \times \bar M$ that maps $(\alpha,\beta)\in M \times \bar M$ to $(\beta, \alpha)$, induces a formal antisymplectic involutive automorphism of  $C^\infty(M \times \bar M, M_{\rm diag})$ and, via $S \otimes T$, an automorphism $I$ of $C^\infty(\T,Z)$. We will give an independent definition of the automorphism $I$ that remains valid in the Poisson case. 
Denote by $E:C^\infty(\T,Z) \to C^\infty(M)$ the composition of the identification mapping from $C^\infty(Z)$ to $C^\infty(M)$ with the restriction to the zero section $Z$ of $\T$. Notice that for any $f \in C^\infty(M) \ E(Sf) =f$ and $E(Tf) = f$. (If $*$ is a star product on a symplectic manifold, $E$ is the unit mapping of the corresponding formal symplectic groupoid.)

\begin{proposition}\label{P:poisscomm}
   Let $F \in C^\infty(\T,Z)$ be a formal function such that $E(F) = 0$. If for any function  $u \in C^\infty(M) \ \{F,Tu\}_\T =0$, then $F = 0$. Similarly, if $\{F,Su\}_\T =0$, then $F = 0$. 
\end{proposition}
\begin{proof} We will prove only the first statement of the Proposition. In local coordinates expand the formal functions $F$ and $Tx^k$ in the series
\[
   F(x,\xi) = \sum_{n=1}^\infty F_n(x,\xi) \text{\quad and \quad } Tx^k = x^k + \sum_{n=1}^\infty T^k_n(x,\xi),
\]
where $F_n$ and $T^k_n$ are the homogeneous components of $F$ and $Tx^k$ of degree $n$ with respect to the fiber variables $\xi$, respectively. Since $\{F,Tx^k\}_\T =0$ for all $k$, we obtain that ${\p F_1}/{\p \xi_k} = 0$ and 
\begin{equation}\label{E:ftn}
       \frac{\p F_n}{\p \xi_k}  = - \sum_{i = 1}^{n-1} \left(\frac{\p F_i}{\p \xi_p}\frac{\p T^k_{n - i}}{\p x^p}  - \frac{\p T^k_{n - i}}{\p \xi_p}\frac{\p F_i}{\p x^p}  \right)
\end{equation}
for all $n \geq 2$. Since $F_1$ is homogeneous of degree 1, we see that $F_1 = 0$. Assume that $F_i = 0$ for all $i <n$. Then it follows from formula (\ref{E:ftn}) that ${\p F_n}/{\p \xi_k} = 0$ for all $k$, whence $F_n = 0$. The statement of the Proposition follows by induction.
\end{proof}
\begin{corollary}
Let $F \in C^\infty(\T,Z)$ be a formal function such that $E(F) = f \in C^\infty(M)$. If for any function  $u \in C^\infty(M) \ \{F,Tu\}_\T =0$, then $F = Sf$. Similarly, if $\{F,Su\}_\T =0$, then $F = Tf$.
\end{corollary}
\begin{proof} We will prove only the first statement of the Corollary.
Since for any functions $u,v \in C^\infty(M)$ the formal functions $Su$ and $Tv$ Poisson commute, the function $F - Sf$ is such that $E(F - Sf) = 0$ and it Poisson commutes with $Tu$  for all $u\in C^\infty(M)$. Therefore, according to Proposition \ref{P:poisscomm}, $F-Sf=0$.
\end{proof}
Let $*$ be a natural star product on a Poisson manifold 
$(M,\{\cdot,\cdot\})$ and $B$ be its automorphism. It follows from Theorem \ref{T:GR} that $B$ can be represented in the form $B = \exp \frac{1}{\nu}X$, where $X$ is a natural differential operator such that $X = 0 \pmod{\nu^2}$.
\begin{proposition}\label{P:symbaut}
For any automorphism $B = \exp \frac{1}{\nu}X$ of a natural star product $*$ on a Poisson manifold $M$ the operator $X$ has a trivial $\sigma$-symbol, $\sigma(X) =0$, and therefore the action of $B$ by conjugation on the algebra $\N$ of natural operators induces the trivial action on the $\sigma$-symbols.
\end{proposition}
\begin{proof}
For $f,g \in C^\infty(M)$ we have that $B(L_f g) = B(f * g) = Bf * Bg = L_{Bf} Bg$, whence $B L_f B^{-1} = L_{Bf}$. On the one hand, $\sigma\left(B L_f B^{-1}\right) = \left(\exp H_{\sigma(X)}\right)\sigma(L_f) = 
\left(\exp H_{\sigma(X)}\right)Sf$.
Since $B = 1 \pmod{\nu}$, we see that $Bf = f \pmod{\nu}$. Therefore, on the other hand,  $\sigma\left(L_{Bf}\right) = \sigma\left(L_f\right) = Sf$. Thus $\left(\exp H_{\sigma(X)}\right)Sf = Sf$ for any $f \in C^\infty(M)$, which implies that $\{\sigma(X),Sf\} = H_{\sigma(X)}Sf = 0$. Since $X = 0 \pmod{\nu^2}$, it follows that $E(\sigma(X))=0$. Using Proposition \ref{P:poisscomm}  we finally obtain that $\sigma(X) = 0$.
\end{proof}
Proposition \ref{P:symbaut} implies the following rigidity property of the equivalence operators of natural deformation quantizations.
\begin{corollary}
If $*$ and $\tilde *$ are equivalent natural star products on $M$ and
$B_1,B_2: (C^\infty(M)[[\nu]],*) \to (C^\infty(M)[[\nu]],\tilde *)$ are equivalence operators of these star products, then the action of the operators $B_1$ and $B_2$ by conjugation on the algebra $\N$ of natural operators induces the same action on the $\sigma$-symbols.
\end{corollary}

Now we will consider the transposition of natural operators on $M$ 
with respect to the pairing (\ref{E:pair}) for some nonvanishing formal density $\rho$ on $M$ and the induced mapping on their $\sigma$-symbols. The mapping $\epsilon: (x,\xi) \mapsto (x, - \xi)$ is a global antisymplectic involutive automorphism of $\T$. It induces the anti-Poisson automorphism $\epsilon^*$ of $C^\infty(\T,Z)$.  For a differential operator $A$ of order $r$ on $M$ its formal transpose $A^t_\rho$ is also of order $r$ and its $r$-symbol is expressed via the $r$-symbol of $A$ as follows, $\s_r(A^t_\rho) = \epsilon^* \s_r(A)$. Therefore, if $A\in \N$ is a natural operator on $M$, then so is $A^t_\rho$, and $\sigma(A^t_\rho)  = \epsilon ^* \sigma(A)$. It follows that the mapping $K_\rho \in \End \D[[\nu]]$ such that
for $A \in \D[[\nu]] \ K_\rho [A] = J_\rho^{-1} A^t_\rho J_\rho$ leaves invariant the algebra $\N$ of natural operators and induces a well-defined mapping on the $\sigma$-symbols that is independent of the choice of the density $\rho$. Represent $J_\rho$ in the form
$J_\rho = \exp \frac{1}{\nu}X$ and denote $I := \exp \left( - H_{\sigma(X)}\right)\circ \epsilon ^*$. Then for $A \in \N$
\begin{equation}\label{E:ik}
 \sigma(K_\rho[A]) = \left(\exp \left( -     H_{\sigma(X)}\right)\circ \epsilon ^*\right)\sigma(A) = I \sigma(A).  
\end{equation}

\begin{theorem}\label{T:last}
The mapping $I$ is an involutive anti-Poisson automorphism of $C^\infty(\T,Z)$ such that $I \circ S =T$ and $I \circ T = S$.
\end{theorem}
\begin{proof} Involutivity of $I$ follows immediately from formula (\ref{E:square}) and Pro\-position \ref{P:symbaut} applied to the modular automorphism $Q_\rho$. Passing to the $\sigma$-symbols on the both sides of (\ref{E:krl}) and using formula (\ref{E:ik}) we obtain that $I \circ T = S$. The identity $I \circ S = T$ follows from the involutivity of $I$.
\end{proof}
Theorem \ref{T:last} has the following obvious corollary.
\begin{corollary}
If $*$ is a natural deformation quantization on a symplectic manifold $M$, then the mapping $I$ is the inverse mapping of the corresponding formal symplectic groupoid.
\end{corollary}

\end{document}